\input amstex
\documentstyle{amsppt}
\document
\magnification=1200
\NoBlackBoxes
\nologo


\centerline{\bf SIXTH PAINLEV\'E EQUATION, UNIVERSAL ELLIPTIC CURVE,}

\smallskip

\centerline{\bf AND MIRROR OF $\bold{P}^2$}

\medskip


\bigskip

\centerline{Yu.I.Manin}

\smallskip

\centerline{\it Max--Planck--Institut f\"ur Mathematik, Bonn, Germany}

\bigskip

\centerline{\bf 0. Introduction}

\medskip

{\bf 0.1. Three approaches to the Painlev\'e equations.}
The differential equations studied in this paper
form a family PVI${}_{\alpha ,\beta , \gamma ,\delta}$
depending on four parameters
$\alpha ,\beta , \gamma ,\delta$, and classically written as:
$$
\frac{d^2X}{dt^2}=\frac{1}{2}\left(
\frac{1}{X}+\frac{1}{X-1}+\frac{1}{X-t}\right)
\left(\frac{dX}{dt}\right)^2 -
\left(
\frac{1}{t}+\frac{1}{t-1}+\frac{1}{X-t}\right)\frac{dX}{dt} +
$$
$$
+\frac{X(X-1)(X-t)}{t^2(t-1)^2}
\left[ \alpha +
\beta\frac{t}{X^2}+\gamma\frac{t-1}{(X-1)^2}+
\delta\frac{t(t-1)}{(X-t)^2}\right]. \eqno{(0.1)}
$$
They were discovered around 1906 and have been approached
from at least three different directions.

\smallskip

{\it a. Study of non--linear ordinary differential equations of the second order
whose solutions have no movable critical points.}

\smallskip

Their classification program was initiated by Painlev\'e,
but he inadvertently omitted (0.1) due to an error
in calculations. It was B. Gambier [G] who completed
Painlev\'e's list and found (0.1).

\smallskip

{\it b. Study of the isomonodromic deformations of linear
differential equations.}

\smallskip

{\it c. Theory of abelian integrals depending on parameters
and taken over chains with algebraic boundary (not necessarily cycles.)}

\smallskip

These two approaches are due to R. Fuchs [F].

\smallskip

In the subsequent development of the theory, relationship
with isomonodromic deformations proved to be most fruitful.
For some recent research and bibliography the reader
may consult [JM], [O1], [H1], [H2].

\smallskip

In this paper I take up the somewhat neglected approach
via abelian integrals and algebraic geometry.

\smallskip

My principal motivation was the desire to understand the
quantum cohomology of $\bold{P}^2$ and to find
an algebraic--geometric object which could be reasonably
called the mirror of $\bold{P}^2$, thus tentatively extending
the scope of the mirror duality discovered for
Calabi--Yau manifolds.

\smallskip

As was explained in a preprint version of [D] (cf.~also 
[DFI] and [H3]), the potential of the quantum cohomology of $\bold{P}^2$
can be reduced by a change of variables  to a particular
solution of the Painlev\'e equation with parameters
$(\alpha ,\beta , \gamma ,\delta )=(\frac{1}{8}, -\frac{1}{8},0,\frac{1}{2}).$
The Painlev\'e transcendents are generally
``new'' functions, but for certain values of parameters all or some solutions
can be expressed through more classical special functions.
Whether this is true for the $\bold{P}^2$--solution referred to above,
seems an open problem. N. Hitchin completely solved
the equation $(\frac{1}{8}, -\frac{1}{8},\frac{1}{8},\frac{3}{8})$ 
in elliptic functions: cf. [H2].

\smallskip

Trying to understand all this, I arrived to the basically
algebraic--geometric picture of all Painlev\'e VI equations,
which in particular suggests that the mirror of $\bold{P}^2$
can be thought of as a pencil of elliptic curves
with labelled sections of order two and an additional,
possibly transcendental, multisection. More precisely,
the Picard--Fuchs equation for the periods  
of the mirror dual Calabi--Yau family is replaced in our
framework by a ``non--homogeneous Picard--Fuchs equation''
satisfied by the Abelian integral from zero to this 
additional multisection (cf. formula (1.5) below.)
It would be important to understand whether this pattern
persists for quantum cohomology of other Fano manifolds.

\smallskip

I will now briefly describe this picture stressing its
geometric aspects. A reader with more analytic
background may prefer to skip the following section.
In the main body of the paper, a prominent role is given to
the uniformization picture using elliptic and modular functions.
In particular, it allows us to reduce the $\bold{P}^2$--equation
to the beautiful form
$$
\frac{d^2z}{d\tau^2}=-\frac{1}{8\pi^2}\wp_z(z,\tau ), \eqno{(0.2)}
$$
where $\wp$ is the Weierstrass function.

\medskip

{\bf 0.2. Algebraic geometry of Painlev\'e VI: a review.}
We will describe a series of constructions which starts with
a pencil of elliptic curves. We work in the category
of complex analytic manifolds, although the most
natural category for this part seems to be that of schemes
over $\roman{Spec}\ \bold{Z}[\frac{1}{2}].$

\smallskip

a. Let $(\pi :E\to B;D_0,\dots ,D_3)$ be a pencil of
compact smooth curves of genus one, with variable absolute invariant,
endowed with four labelled sections $D_i$ such that
if any one of them is taken as zero, the others
will be of order two.

\smallskip

We will call $E$ {\it a configuration space} of PVI
(common for all values of parameters.) Solutions of
all equations will be represented by some multisections
of $\pi .$

\smallskip

b. Let $\Cal{F}$ be the subsheaf of the sheaf of
vertical 1-forms $\Omega^1_{E/B}(D_3)$ on $E$
with pole at $D_3$ and residue 1 at this pole.
It is an affine twisted version of $\Omega^1_{E/B}$
which is the sheaf of sections of the relative
cotangent bundle $T^*_{E/B}.$ Similarly,
$\Cal{F}$ itself ``is'' the sheaf of sections
of an affine line bundle $F=F_{E/B}$ on $E.$
More precisely, we construct such a bundle
$\lambda :F\to E$ and a form $\nu_F\in\Gamma (F, \Omega^1_{F/B}(\lambda^{-1}(D_3))$
such that the map
$$
\{\roman{local\ section}\ s\ \roman{of}\ F\}\mapsto
s^*(\nu_F)
$$
identifies the sheaf of sections of $F/E$ with $\Cal{F}.$

\smallskip

We will call $F$ {\it a phase space} for PVI (again, common
for all parameter values.)

\smallskip

c. $E$ carries a distinguished family of algebraic curves
transversal to the fibers of $E$: considered as multisections of
$E/B$ they are of finite order (if any of $D_i$ is chosen as zero.)
It is important that each curve of this family has a canonical
lifting to $F$ (for its description, see the main text, formulas
(2.12) and (2.29).)

\smallskip

d. $F$ carries a closed 2-form $\omega^{(0)}$ which can be characterized
by the following two properties:

\smallskip

{\it i). The vertical part of $\omega^{(0)}$, i.~e. its restriction
to the fibers of $\pi\circ\lambda : F\to B$, coincides
with $d_{F/B}(\nu_F).$}

\smallskip

{\it ii). Any canonical lift to $F$ of a connected multisection
of finite order
of $E\to B$, referred to above, is a leaf of the
null--foliation of $\omega^{(0)}$.}

\smallskip

e. $E$ also carries four distinguished closed two--forms
$\omega_0,\dots ,\omega_3.$ They are determined,
up to multiplication by a constant, by the following
properties.

\smallskip

{\it iii). The divisor of $\omega_i$ is $\dfrac{D_jD_kD_l}{D_i^3}$
where $\{ i,j,k,l\}=\{ 0,1,2,3\}.$}

\smallskip

{\it iv). Identify the sheaves $\Omega^2_E$ and $\pi^*(\Omega^1_{E/B})^{\otimes 3}$
on $E$ using the Kodaira--Spencer isomorphism
$\pi^*(\Omega^1_B)\cong (\Omega^1_{E/B})^{\otimes 2}$
and the exact sequence $0\to \pi^*(\Omega^1_B)\to
\Omega^1_E \to \Omega^1_{E/B} \to 0.$  Then the image of
$\omega_i$ in $\pi^*(\Omega^1_{E/B})^{\otimes 3}$
considered in the formal neighborhood of $D_i$ is the cube
of a vertical 1--form with a constant residue along $D_i.$}

\smallskip

The affine space $P_0:= \omega^{(0)}+\sum_{i=0}^3 \bold{C}\lambda^*(\omega_i)$
of closed two--forms on $F$ is our version of the
{\it moduli space of the PVI equations} replacing the
classical $(\alpha ,\beta ,\gamma ,\delta )$--space.

\smallskip

We can now summarize our definition of PVI equations and
their solutions.

\smallskip

\proclaim{\quad 0.2.1. Definition} a). A Painlev\'e two--form
on $F$ is a point $\omega\in P_0.$

b). The Painlev\'e foliation corresponding to $\omega$
is the null--foliation of $\omega$.

c). The solutions to the respective Painlev\'e equation
are the leaves of this foliation (in the Hamiltonian
description), or their projections on $E$. 

The form $\omega^{(0)}$ corresponds to
$(\alpha ,\beta ,\gamma ,\delta )=(0,0,0,\frac{1}{2}).$
\endproclaim

To obtain (0.1), we must specialize this description
to the (projectivized) family
$$
E_t:\ Y^2=X(X-1)(X-t) \mapsto t\in B:=\bold{P}^1\setminus\{0,1,\infty\}
\eqno{(0.3)}
$$
and look at the variation of $X$ along solutions.

\medskip

{\bf 0.3. Plan of the paper.} In \S 1 we reproduce
R. Fuchs's description of (0.1) in terms of elliptic
integrals and deduce from it an
analytic form of the Painlev\'e equations
involving Weierstrass $\wp$--punction. As an application,
we derive the elementary symmetries of PVI and
introduce the Landin transform.

The key \S 2 is devoted to the Hamiltonian structure
of PVI and contains proofs of all claims made in 0.2 above.

In \S 3 we establish the relationship of our Hamiltonian picture
with that of Okamoto [O2] and sketch
Okamoto's treatment of the hidden $W(D_4)$--symmetry of PVI.
We also review some known solutions.

This symmetry nicely explains, why Hitchin was able to
solve his 
$(\frac{1}{8}, -\frac{1}{8},\frac{1}{8},\frac{3}{8})$--equation.

\medskip

{\bf 0.4. Further plans.} The geometric setting advocated in this
paper furnishes a convenient framework for the treatment
of the following subject matters:

a. Three--dimensional Frobenius manifolds, including the
quantum cohomology of $\bold{P}^2.$

b. Geometry of the degenerations of PVI to PV, \dots, PI.

c. Generalizations to higher genus and isomonodromic
deformations with many singular points.

I hope to return to these problems in future publications.

\medskip

{\it Acknowledgement.} I am very grateful to Andrey Levin for
consultations and help with elliptic functions.

\vskip 1cm

\centerline{\bf \S 1. PVI and elliptic functions}

\bigskip

\proclaim{\quad 1.1. Theorem (R. Fuchs, 1907)} The equation (0.1)
can be written in the form
$$
t(1-t)\left[t(1-t)\frac{d^2}{dt^2}+(1-2t)\frac{d}{dt}-\frac{1}{4}\right]
\int_{\infty}^{(X,Y)}\frac{dx}{\sqrt{x(x-1)(x-t)}}=
$$
$$
=\alpha Y+\beta\frac{tY}{X^2}+\gamma\frac{(t-1)Y}{(X-1)^2}
+(\delta -\frac{1}{2})\frac{t(t-1)Y}{(X-t)^2}  \eqno{(1.1)}
$$
where $Y^2=X(X-1)(X-t).$
\endproclaim

{\bf Proof.} First, let us clarify the meaning of (1.1).
Consider the family of alliptic curves $E\to B$ parametrized
by $t\in \bold{P}^1\setminus \{0,1,\infty\} := B:$ the curve
$E_t$ is the projective closure of $Y^2=X(X-1)(X-t).$
Points at infinity of $\{E_t\}$ form a section $D_0$ of this family which
is the zero section for the standard group law on fibers.
Choose in $E_t(\bold{C})$ a path from $D_0(t)$ to the point $(X(t),Y(t))$
of a local section. The operator
$$
L_t:=t(1-t)\frac{d^2}{dt^2}+(1-2t)\frac{d}{dt}-\frac{1}{4} \eqno{(1.2)}
$$
annihilates the periods $\dsize\int \frac{dx}{y}$ along closed paths in
$E_t(\bold{C})$ because
$$
\left[t(1-t)\frac{\partial^2}{\partial t^2}+(1-2t)\frac{\partial}{\partial t}-\frac{1}{4}\right]\frac{d_{E/B}x}{y}=
\frac{1}{2}d_{E/B}\frac{y}{(x-t)^2}
\eqno{(1.3)}
$$
where we put $\dfrac{\partial}{\partial t}(x)=0$ and $d_{E/B}t=0.$
Applying $L_t$ to $\dsize\int_{\infty}^{(X,Y)} \frac{dx}{y}$ we get
$\left.\dfrac{1}{2}\dfrac{y}{(x-t)^2}\right|^{(X,Y)}_{\infty}$
plus the contribution of the boundary sections
which together with the right hand side of (1.1)
amounts to (0.1).

\medskip

{\bf 1.2. $\mu$--equations.} The equation (1.1) is an instance
of a general construction which was used in [M] to prove
the functional Mordell conjecture. We will recall it now.

\smallskip

A {\it $\mu$--equation} is a system of non--linear PDE
in which independent variables are (local) coordinates
on a manifold $B$ and unknown functions are represented
by a section $s$ of a family of abelian varieties
(or complex tori) $\pi : A\to B.$ To write this system explicitly,
assume $B$ small enough so that $\pi_*(\Omega^1_{A/B})$
and $\Cal{D}_B$ (sheaf of differential operators on $B$)
are $\Cal{O}_B$--free, and make the following choices:

\smallskip

a. An $\Cal{O}_B$--basis of vertical 1--forms
$\omega_1,\dots ,\omega_n\in \Gamma (B, \pi_*(\Omega^1_{A/B})).$

\smallskip

b. A system of generators of the $\Cal{D}_B$--module of the Picard--Fuchs
equations
$$
\sum_{i=1}^n L_i^{(j)}\int_{\gamma}\omega_i=0,\quad j=1,\dots ,N,
\eqno{(1.4)}
$$
where $\gamma$ runs over families of closed paths in the fibers
spanning $H_1(B_t)$.

c. A family of meromorphic functions $\Phi^{(j)},\ j=1,\dots ,N$
on $A.$

\smallskip

The respective $\mu$--equation for a local (multi)--section
$s:\ B\to A$ reads then
$$
\sum_{i=1}^n L_i^{(j)}\int_{0}^{s}\omega_i= s^*(\Phi^{(j)}),\quad j=1,\dots ,N,
\eqno{(1.5)}
$$
where $0$ denotes the zero section.

\smallskip

One drawback of (1.5) is its dependence on arbitrary choices.
Clearly, this can be reduced by taking account of the
transformation rules with respect to the changes of 
various generators. For elliptic pencils, the result
takes a very neat form.

\medskip

{\bf 1.3. Elliptic $\mu$--equations.} Let again $E\to B$
be a non--constant one--dimensional family of elliptic curves.
We temporarily keep the assumption that $\pi_*(\Omega^1_{E/B})$
and the tangent sheaf $\Cal{T_B}$ are free.
For any symbol of order two $\sigma\in S^2(\Cal{T_B})$
and any generator $\omega$ of $\pi_*(\Omega^1_{E/B})$
denote by $L_{\sigma ,\omega}$ the Picard--Fuchs operator on $B$
with the symbol $\sigma$ annihilating all periods of $\omega .$

\smallskip

\proclaim{\quad 1.3.1. Lemma} For any local section $s$,
the expression $\dsize L_{\sigma ,\omega}\int_{0}^s\omega$
is $\Cal{O}_B$--bilinear in $\sigma$ and $\omega .$
\endproclaim

{\bf Proof.} Obviously,
$$
L_{f\sigma ,\omega}=fL_{\sigma ,\omega},\quad
L_{\sigma ,g\omega}=gL_{\sigma ,\omega}\circ g^{-1},
$$
where $f,g$ are functions on $B$. The lemma follows.

\smallskip

Thus the expression
$$
\mu (s):= \left( L_{\sigma ,\omega}\int_{0}^s\omega\right)
\otimes\sigma^{-1}\otimes{\omega}^{-1}\in
S^2(\Omega^1_B)\otimes(\pi_*\Omega^1_{E/B})^{-1} \eqno{(1.6)}
$$
depends only on $s$ and is compatible with restrictions to open subsets
of $B$. This means that the natural domain of the
right hand sides for elliptic $\mu$--equations
is the set of meromorphic sections $\Phi$ of the sheaf
$\pi^*\left[S^2(\Omega^1_B)\otimes(\pi_*\Omega^1_{E/B})^{-1}\right].$

\smallskip

Notice that the Kodaira--Spencer isomorphism (and eventually a choice
of the theta--characteristic of $B$) allows us to identify
$\Phi$ with a meromorphic section of
$(\Omega^1_{E/B})^3$ or $\pi^*(\Omega^1_B)^{3/2}$ as well.

\smallskip

We will now lift the Fuchs--Painlev\'e equation (1.1)
to the classical covering space, which in particular will make transparent
the nature of its right hand side.

\medskip

{\bf 1.3.2. Uniformization.} Consider the family of elliptic
curves parametrized by the upper half--plane $H$:
$E_{\tau}:= \bold{C}/(\bold{Z}+\bold{Z}\tau )\mapsto
\tau\in H$. Recall that
$$
\wp (z,\tau ):= \frac{1}{z^2}+\sum{}^{\prime}
\left(\frac{1}{(z+m\tau +n)^2}-\frac{1}{(m\tau +n)^2}\right) ,      
\eqno{(1.7)}
$$
$$
\wp_z (z,\tau )= 
-2\sum\frac{1}{(z+m\tau +n)^3} . \eqno{(1.8)}
$$
We have
$$
\wp_z (z,\tau )^2=4(\wp (z,\tau )-e_1(\tau ))(\wp (z,\tau)-e_2(\tau ))
(\wp (z,\tau )-e_3(\tau ))
\eqno{(1.9)}
$$
where
$$
e_i(\tau )=\wp (\frac{T_i}{2},\tau ),\ (T_0,\dots ,T_3)=
(0,1,\tau ,1+\tau )
 \eqno{(1.10)}
$$
and $e_1+e_2+e_3=0.$ Functions $\wp$ and $\wp_z$
are invariant with respect to the shifts $\bold{Z}^2:
(z,\tau )\mapsto (z+m\tau +n,\tau)$
and behave in the following way under the full modular
group $\Gamma :$
$$
\wp\left(\frac{z}{c\tau +d},\frac{a\tau +b}{c\tau +d}\right)=
(ct+d)^2\wp (z,\tau ), \eqno{(1.11)}
$$
$$
\wp_z\left(\frac{z}{c\tau +d},\frac{a\tau +b}{c\tau +d}\right)=
(ct+d)^3\wp_z (z,\tau ). \eqno{(1.12)}
$$
Consider now the morphism of families $\varphi :\ \{E_{\tau}\}\to \{E_t\}$
induced by
$$
(z,\tau )\mapsto \left( X=\frac{\wp (z,\tau)-e_1}{e_2-e_1},
Y=\frac{\wp_z (z,\tau)}{2(e_2-e_1)^{3/2}},
t=\frac{e_3-e_1}{e_2-e_1}\right) .
\eqno{(1.13)}
$$
This is a Galois covering with the group $\Gamma (2)\ltimes\bold{Z}^2.$
We have
$$
\varphi^*\left(\frac{d_{E/B}X}{Y} \right) =
2(e_2-e_1)^{1/2}d_{E/H}z.  \eqno{(1.14)}
$$
In the future formulas of this type we will omit $\varphi^*$
and denote differentials over a base $B$ by $d_{\downarrow}.$
For instance,
$d_{\downarrow}\left(\dfrac{z}{c\tau +d} \right) =
\dfrac{d_{\downarrow}z}{c\tau +d}$, whereas
$d\left(\dfrac{z}{c\tau +d} \right) =
\dfrac{dz}{c\tau +d}-\dfrac{czd\tau}{(c\tau +d)^2}.$

\smallskip

It follows from (1.14) that if we denote by $\gamma_1$
(resp. $\gamma_2$) the image of [0,1] (resp. [0,1]$\tau$)
in $\{E_t\}$, then
$$
\int_{\gamma_1}\frac{d_{\downarrow}X}{Y}=2(e_2-e_1)^{1/2},\quad
\int_{\gamma_2}\frac{d_{\downarrow}X}{Y}=2\tau (e_2-e_1)^{1/2}
\eqno{(1.15)}
$$
so that the operator $L_t$ from (1.2) annihilates
periods (1.15) as functions of $\tau$.

\smallskip

\proclaim{\quad 1.4. Theorem} A lift of (1.1) to the
$(z,\tau )$--space $\bold{C}\times H$ reads:
$$
\frac{d^2z}{d\tau^2}=\frac{1}{(2\pi i)^2}\sum_{j=0}^3
\alpha_i\wp_z(z+\frac{T_i}{2},\tau )
\eqno{(1.16)}
$$
where
$$
(\alpha_0,\dots ,\alpha_3):=(\alpha ,-\beta ,\gamma ,
\frac{1}{2}-\delta ).
\eqno{(1.17)}
$$
\endproclaim
 
{\bf Proof.} Following the lead of no. 1.3, we will
directly calculate the $\mu$--equation for $\{E_{\tau}\}$,
choosing $\omega = d_{\downarrow}z$ (instead of $d_{\downarrow}X/Y$)
and $\sigma =\dfrac{d^2}{d\tau^2}$ (instead of $t^2(1-t)^2
\dfrac{d^2}{dt^2}$.) Since periods of $ d_{\downarrow}z$
are generated by 1 and $\tau$, the relevant Picard--Fuchs operator
is simply $\dfrac{d^2}{d\tau^2}$. From the Lemma 1.3.1 and
(1.15) it follows that
$$
t(1-t)L_t\circ 2(e_2-e_1)^{1/2}=Z(\tau )\frac{d^2}{d\tau^2}.
$$
Using (1.13) and comparing symbols, we see that
$$
Z(\tau )=2\left(\frac{e_3-e_1}{e_2-e_1} \right)^2
\left(\frac{e_3-e_2}{e_2-e_1} \right)^2
\frac{(e_2-e_1)^4}{9(e_1e_2^{\prime}-e_2e_1^{\prime} )^2}
(e_2-e_1)^{1/2}=
$$
$$
=\frac{2}{9}\,\frac{\prod_{i>j}(e_i-e_j)^2}{(e_1e_2^{\prime}-e_2e_1^{\prime} )^2}(e_2-e_1)^{-3/2}.
\eqno{(1.18)}
$$
Since $e_1+e_2+e_3=0,$ we can replace $(e_1e_2^{\prime}-e_2e_1^{\prime} )^2$
by $(e_ie_j^{\prime}-e_je_i^{\prime} )^2$ for any $i\neq j.$
It follows that
$$
C:=\frac{\prod_{i>j}(e_i-e_j)^2}{(e_1e_2^{\prime}-e_2e_1^{\prime} )^2}
$$
is a modular function for the full modular group without
zeroes and poles, hence a constant. A calculation with
theta--functions, here omitted, for which I am grateful
to A. Levin, shows that $C=-9\pi^2,$ so that finally
$$
t(1-t)L_t\int_{\infty}^{(X(t),Y(t))}\frac{d_{\downarrow}x}{y}=
-2\pi^2(e_2-e_1)^{-3/2}\frac{d^2}{d\tau^2}
\int_0^{z(\tau )}d_{\downarrow}z
\eqno{(1.19)}
$$
for the respective sections. We can now consecutively
compare the summands in the right hand side of (1.1)
with those in (1.16). The first summand
gives
$$
\alpha Y=\frac{\alpha}{2}(e_2-e_1)^{-3/2}\wp_z (z,\tau ).
$$
For the remaining ones we have to use the addition formulas
$$
\wp_z(z+\frac{T_i}{2},\tau )=
-\frac{(e_i-e_j)(e_i-e_k)}{(\wp (z,\tau )-e_i)^2}
\wp_z (z,\tau ),\quad \{i,j,k\} = \{1,2,3\},
$$
so that, say, for $i=3$ we get
$$
(\delta -\frac{1}{2})\frac{t(t-1)Y}{(X-t)^2}=
(\delta -\frac{1}{2})\frac{(e_3-e_1)(e_3-e_2)}{(e_2-e_1)^2}\cdot
\frac{\wp_z (z,\tau )}{2(e_2-e_1)^{3/2}}\cdot
\frac{(e_2-e_1)^2}{(\wp (z,\tau )-e_3)^2}=
$$
$$
=-\frac{1}{2}\,(\delta -\frac{1}{2})\,(e_2-e_1)^{-3/2}\cdot
\frac{-(e_3-e_1)(e_3-e_2)}{(\wp (z,\tau )-e_3)^2}\wp_z (z,\tau )=
$$
$$
=-\frac{1}{2}(\delta -\frac{1}{2})(e_2-e_1)^{-3/2}
\wp_z(z+\frac{1+\tau}{2},\tau ).
$$
The remaining two summands are treated similarly.
This finishes the proof.

\smallskip

As the first application, we can now describe the space
of the right hand sides of Painlev\'e--Fuchs equations
in a model--independent way (compare Introduction, 0.2e, iv)).

\medskip

{\bf 1.5. PVI on an arbitrary elliptic pencil.} We put ourselves in the
setting of 0.2a. As was explained in 1.3.1, for an
invariant $\mu$--equation (as (1.6)) the right hand side can be
considered as a meromorphic section of $\pi^*(S^3(\Omega^1_{E/B}))$.
The space of the right hand sides of (1.16) written in the
invariant form is generated by four cubic differentials
$\wp_z(z+\frac{T_i}{2},\tau )(d_{\downarrow}z)^3$.
Looking at their Laurent series near $z+\dfrac{T_i}{2}=0$ one easily sees
that they are cubes of a formal differential with a constant residue
along $D_i:=\dfrac{T_i}{2}\ \roman{mod} (1,\tau )$,
and that this property together with identification of their
divisors as $\dfrac{D_jD_kD_l}{D_i^3}$ characterizes them
up to a multiplicative constant.

\smallskip

In the Theorem 2.5 below we will give a Hamiltonian interpretation
of this space.

\medskip

{\bf 1.6. $S_4$--symmetry and the Landin transform.} As the first application
of 1.5 and (1.16) we will construct some natural
transformations of PVI. For much deeper hidden symmetries,
see \S 3.

\smallskip

a. {\it The classical $S_4$--symmetry.} Isomorphisms of $(E,D_i)$
which do not conserve the labelling of $D_i$ induce transformations
of PVI permuting $\alpha_i$. In the form (1.16), they act on
solutions as compositions of the transformations of two types:
$\dsize (z,\tau )\mapsto \left( \frac{z}{cz+\tau},\frac{a\tau +b}{c\tau +d}
\right)$ indexed by cosets $\Gamma /\Gamma (2)$, and
$\dsize (z,\tau )\mapsto (z+\frac{T_i}{2},\tau )$ shifting the
zero section.

\smallskip

b. {\it The Landin transform.} From (1.8) one easily deduces
Landin's identity
$$
\wp_z(z,\frac{\tau}{2})=
-2\left[\sum\frac{1}{(z+2m\frac{\tau}{2} +n)^3}+   
\sum\frac{1}{(z+\frac{\tau}{2}+2m\frac{\tau}{2} +n)^3} \right]=
$$
$$
=\wp_z(z,\tau )+\wp_z(z+\frac{\tau}{2},\tau ).
$$
Hence if $z(\tau )$ is a solution to PVI with parameters
$(\alpha_0,\alpha_1,\alpha_0,\alpha_1)$, we have
$$
\frac{d^2z(\tau )}{d\tau^2}=
\alpha_0 [\wp_z(z,\tau )+\wp_z(z+\frac{\tau}{2},\tau )]+
\alpha_1 [\wp_z(z+\frac{1}{2},\tau )+\wp_z(z+\frac{1+\tau}{2},\tau )]=
$$
$$
=\frac{1}{4}\,\frac{d^2z(\tau )}{d(\tau /2)^2}
=\alpha_0 \wp_z(z,\frac{\tau}{2})+\alpha_1 \wp_z(z+\frac{1}{2},\frac{\tau}{2}),
$$
that is, $z(2\tau )$ is a solution to PVI with parameters
$(4\alpha_0,4\alpha_1,0,0).$ The converse statement is true as well.
In this way we get the following bijections between the sets
of solutions to (1.16):
$$
(\alpha_0,\alpha_1,\alpha_0,\alpha_1)\leftrightarrow
(4\alpha_0,4\alpha_1,0,0) 
\eqno{(1.19)}
$$
and in particular
$$
(\alpha_0,0,\alpha_0,0) \leftrightarrow
(4\alpha_0,0,0,0).
\eqno{(1.20)}
$$
Of course, we can combine these correspondences with permutations.
In this way, Hitchin's equation reduces in two steps to
$$
\frac{d^2z}{d\tau^2}=-\frac{1}{2\pi^2}\wp_z(z,\tau ),
\eqno{(1.21)}
$$
whereas the $\bold{P}^2$--equation reduces to (0.2).

\medskip

{\bf 1.7. Remark.} A straightforward generalization of (1.16)
is the following infinite--dimensional family
of $\mu$--equations:
$$
\frac{d^2z}{d\tau^2}=\frac{1}{(2\pi i)^2}
\sum_{\zeta}\alpha_{\zeta}\wp_z(z+\zeta ,\tau ),
\eqno{(1.22)}
$$
where $\zeta$ runs over representatives of
$(\bold{Q}+\bold{Q}\tau )/(\bold{Z}+\bold{Z}\tau )$,
and $\alpha_{\zeta}=0$ for almost all $\zeta$.
Most of the results of this paper readily
extend to (1.22)

\vskip 1cm

\centerline{\bf \S 2. Hamiltonian structure}

\bigskip

{\bf 2.1. The time--dependent Hamiltonian.} The PVI--equation written
as in (1.16) has an obvious time--dependent Hamiltonian form:
$$
\frac{dz}{d\tau}=\frac{\partial\Cal{H}}{\partial y},\
\frac{dy}{d\tau}=-\frac{\partial\Cal{H}}{\partial z},
\eqno{(2.1)}
$$
where
$$
\Cal{H}:=\frac{y^2}{2}-\frac{1}{(2\pi i)^2}
\sum_{j=0}^{3}\alpha_j\wp (z+\frac{T_j}{2},\tau ).
\eqno{(2.2)}
$$
To understand the geometric meaning of these equations, we will
extend the action of $\Gamma (2)\ltimes\bold{Z}^2$ to the
$(y,z,\tau )$--space in a way compatible with (2.1), (2.2).
We start with recalling the general Hamiltonian formalism.

\medskip

{\bf 2.2. Hamiltonian formalism.} {\it a. Non--degenerate case.}
Let $X$ be a manifold, $\pi\in\Gamma (X,\wedge^2\Cal{T}_X),\,
\omega\in\Gamma (X,\wedge^2\Cal{T}^*_X)$. The natural
integrability conditions for such tensors are

\smallskip

for $\pi$:\quad $\{f,g\}_{\pi}:=\pi (df,dg)\ 
\roman{satisfies\ the\ Jacobi\ identity}$;

for $\omega$:\quad $d\omega =0.$

\smallskip

If both $\pi$ and $\omega$ are nowhere degenerate, we can write
the {\it compatibility condition} for them meaning
that they define mutually inverse isomorphisms
$\Cal{T}_x\overset \tilde{\omega}\to{
\underset \tilde{\pi}\to\rightleftarrows}\Cal{T}^*_X$.
This relation is a bijection compatible with the two
integrability conditions, which establishes the equivalence between
the non--degenerate Poisson structures $\pi$ on $X$
and the symplectic structures $\omega$ on $X$, so that we can
write the relevant Poisson bracket as $\{f,g\}_{\omega}$ 
as well. Any function $\Cal{H}$ on $X$ (time--independent
Hamiltonian) defines a flow on $X$ endowed with $\pi$ or
$\omega$. This flow has respectively two equivalent
descriptions:

\smallskip

Poisson:\quad $\dfrac{df}{dt}=\{\Cal{H},f\}_{\omega}$, $f$
being any function on $X$;

symplectic: graphs of the flow lines in the extended
phase space $X\times A^1_t$ are leaves of the null--foliation
of the closed form $\roman{pr}^*_X\omega -d\Cal{H}\wedge dt.$

\medskip

{\it b. Degenerate case.} Here the two structures diverge,
and the natural compatibility relation ceases to
be a bijection.

\smallskip

A tensor $\pi\in\Gamma (X,\wedge^2\Cal{T}_X)$ of constant rank
defines the subbundle $\roman{Ker}\,\tilde{\pi}\subset\Cal{T}^*_X$
and the orthogonal distribution $(\roman{Ker}\,\tilde{\pi})^{\perp}
\subset\Cal{T}_X$.
If in addition $\pi$ is Poisson, then

\smallskip

i). $(\roman{Ker}\,\tilde{\pi})^{\perp}$ is integrable,
i.~e. it defines a foliation called the symplectic foliation of $\pi$.

ii). On the leaves of this foliation, $\pi$ induces a nondegenerate
Poisson, or equivalently, symplectic structure.

\smallskip

On the other hand, a tensor $\omega\in\Gamma (X,\wedge^2\Cal{T}_X^*)$
of constant rank directly defines the distribution
$\roman{Ker}\,\tilde{\omega}\subset\Cal{T}_X$, and if $\omega$
is closed, then

\smallskip

i${}^\prime$). $\roman{Ker}\,\tilde{\omega}$ is integrable; its leaves form
the null--foliation of $\omega$.

ii${}^\prime$). $\omega$ induces a symplectic structure on the leaves of any
foliation transversal to the null--foliation of
$\omega$ and having the complementary dimension.

\smallskip

We will now call $\pi$ and $\omega$ {\it compatible}, if
$\Cal{T}_X=(\roman{Ker}\,\tilde{\pi})^{\perp}\oplus
\roman{Ker}\,\tilde{\omega}$, and if in addition,
$\pi$ and $\omega$ induce the same symplectic structure
on the symplectic leaves of $\pi$.

\smallskip

In the remaining part of this paper, we will be interested
only in the (degenerate) symplectic picture $(X,\omega )$
considered as a generalization of the extended
phase space. The leaves of the null--foliation
will be for us solutions to a Hamiltonian
system. The following simple Proposition shows that
a particular case of this picture encodes the
classical formalism of Hamiltonian equations
with many times and time--dependent Hamiltonians.

\smallskip

\proclaim{\quad 2.2.1. Proposition} Let $X=X_0\times B,\
(p_i,q_i), i=1,\dots ,n,$ be coordinates on $X_0$,
$(t_1,\dots ,t_m)$ coordinates on $B$. Let 
$\omega_0=\sum_{i=1}^ndp_i\wedge dq_i$ be a non--degenerate
symplectic form on $X_0$, and
$\omega =\sum_{i=1}^ndp_i\wedge dq_i-\sum_{j=1}^m d\Cal{H}_j\wedge dt_j$
be a closed form of the constant rank $2n$, where
$\Cal{H}_j=\Cal{H}_j(p,q,t)$ are functions on $X$. Then we have:

\smallskip

a). Leaves of the null--foliation of $\omega$ form an \'etale covering
of $B$ iff the Hamiltonians $\Cal{H}_j$ satisfy the integrability
condition
$$
\forall j,k,\quad \{\Cal{H}_j,\Cal{H}_k\}_{\omega_0}=
\partial_{t_j}\Cal{H}_k - \partial_{t_k}\Cal{H}_j
$$
(empty for $m=1$), where the Poisson bracket is taken at constant
times.

\smallskip

b). The equations of motion expressing variation of $p_i,q_i$
along the leaves are
$$
\frac{\partial p_i}{\partial t_j}=-\frac{\partial\Cal{H}_j}{\partial q_i},\
\frac{\partial q_i}{\partial t_j}=\frac{\partial\Cal{H}_j}{\partial p_i}.
$$

\endproclaim

{\bf Proof.} Leaves of the null-foliation form an \'etale
covering of $B$, iff the null--distribution is spanned
by lifts of the  basic vector fields $\partial_{t_j}$:
$$
\Cal{D}_j=\sum_{i=1}^n A_i^{(j)}\partial_{p_i} +
          \sum_{i=1}^n B_i^{(j)}\partial_{q_i} +\partial_{t_j},\
j=1,\dots , m.
$$
We have then:

$\omega (\Cal{D}_j,\partial_{p_i})=0\ \Leftrightarrow \
B_i^{(j)}=\dfrac{\partial\Cal{H}_j}{\partial p_i}$,

\smallskip

$\omega (\Cal{D}_j,\partial_{q_i})=0\ \Leftrightarrow \
A_i^{(j)}=-\dfrac{\partial\Cal{H}_j}{\partial q_i}$,

\smallskip

$\omega (\Cal{D}_j,\partial_{t_k})=0\ \Leftrightarrow \
-\sum_{i=1}^n A_i^{(j)}\partial_{p_i}\Cal{H}_k -
          \sum_{i=1}^n B_i^{(j)}\partial_{q_i}\Cal{H}_k -
\partial_{t_j}\Cal{H}_k + \partial_{t_k}\Cal{H}_j =0.$

\smallskip

Finally, for any observable $f$, the equations of motion are
$$
\frac{\partial }{\partial t_j}(f|_L)=(\Cal{D}_jf)|_L
$$
where $L$ is any leaf of the foliation.

\medskip

Proposition 2.2.1 generally furnishes a too simplified picture. Not only Hamiltonians but the constant time slices, together
with their symplectic structure, might become time--dependent
(especially in the analytic context). Still worse, projection onto a
time manifold $B$ may not be a part of the data. Even
a specific transversal foliation need not be present.

\medskip

{\bf 2.3. PVI revisited.} Looking at (2.1) and (2.2),
we see that in the $(y,z,\tau )$--space solutions of a
particular PVI form the null--foliation of
$$
\omega =\omega(\alpha_0,\dots ,\alpha_3):=
2\pi i(dy\wedge dz-d\Cal{H}\wedge d\tau )=
$$
$$
=2\pi i(dy\wedge dz-ydy\wedge d\tau ) +
\frac{1}{2\pi i}
\sum_{j=0}^{3}\alpha_j\wp_z (z+\frac{T_j}{2},\tau )dz\wedge d\tau .
\eqno{(2.3)}
$$
(The extra factor $2\pi i$ makes $\omega$ defined over $\bold{Q}$, cf. below.)

\smallskip

\proclaim{\quad 2.3.1. Proposition} The standard action of
$\Gamma (2)$ (resp. $\bold{Z}^2$) on $\bold{C}\times H$
has a unique extension to $\bold{C}\times \bold{C}\times H$
leaving $(2.3)$ invariant:
$$
(y,z,\tau )\mapsto \left( y(c\tau +d)-cz,\frac{z}{c\tau+d},
\frac{a\tau +b}{c\tau +d}\right),
\eqno{(2.4)}
$$
$$
(y,z,\tau )\mapsto (y+m, z+m\tau +n,\tau ).
\eqno{(2.5)}
$$
\endproclaim

{\bf Proof.} Let $(y,z,\tau )\mapsto \left( \tilde{y},
\tilde{z}=\dfrac{z}{c\tau+d},
\tilde{\tau}=\dfrac{a\tau +b}{c\tau +d}\right),$
$d\tilde{y}=Ady+Bdz+Cd\tau ,$ be a transformation from $\Gamma (2)$
preserving the form of (2.3):
$$
\omega\equiv
2\pi i(d\tilde{y}\wedge d\tilde{z}-\tilde{y}d\tilde{y}\wedge d\tilde{\tau}) +
\frac{1}{2\pi i}
\sum_{j=0}^{3}\alpha_j\wp_z (\tilde{z}+\frac{T_j}{2},\tilde{\tau} )
d\tilde{z}\wedge d\tilde{\tau }.
\eqno{(2.6)}
$$
From (1.12) it follows that the terms in (2.6) involving the Weierstrass
function are automatically invariant. Comparing coefficients
of $dy\wedge dz$ at both sides of (2.6), one sees that
$A=c\tau +d.$ Comparing coefficients of $dy\wedge d\tau,$ one
then finds $\tilde{y}=y(c\tau +d)-cz,$ which gives $B=-c,C=cy.$
Finally, one checks the vanishing of the relevant part of
the coefficient of $dz\wedge d\tau .$ This proves (2.4);
(2.5) is checked similarly.

\medskip

We will now construct a function of $z,\tau$ behaving in the
same way as $y$ in (2.4), (2.5).

\smallskip

Namely, consider the theta--function
$$
\theta (z,\tau )=\sum_{n\in\bold{Z}}\roman{exp}\,(\pi in^2\tau + 
2\pi inz).
$$
It has zeroes of the first order at $z\equiv\dfrac{1+\tau}{2}\,\roman{mod}\,(1,\tau )$
and satisfies the following functional equations under
the action of $\Gamma (2)$ and $\bold{Z}^2$:
$$
\theta \left( \frac{z}{c\tau+d},
\frac{a\tau +b}{c\tau +d}\right) =
\zeta (c\tau +d)^{1/2}\roman{exp}\,\left( \pi ic\,\frac{z^2}{c\tau +d}\right)
\theta (z,\tau),
\eqno{(2.7)}
$$
$$
\theta (z+m\tau +n)=\roman{exp}\,(-\pi im^2\tau -2\pi imz)\,\theta (z,\tau ),
\eqno{(2.8)}
$$
where $\zeta$ is a root of unity of degree eight.
Therefore the function $v(z,\tau):=-\dfrac{1}{2\pi i}\dfrac{\theta_z}{\theta}
(z,\tau )$ has poles of the first order with residue $-\dfrac{1}{2\pi i}$
at $z\equiv\dfrac{1+\tau}{2}\,\roman{mod}\,(1,\tau )$ and satisfies
$$
v\left( \frac{z}{c\tau+d},
\frac{a\tau +b}{c\tau +d}\right) =v(z,\tau)(c\tau +d)-cz,
\eqno{(2.9)}
$$
$$
v(z+m\tau +n)=v(z,\tau )+m.
\eqno{(2.10)}
$$
Comparing this to (2.4), we find finally:

\smallskip

\proclaim{\quad 2.3.2. Proposition} The vertical (over $H$)
differential
$$
\nu :=\left(2\pi iy+\frac{\theta_z}{\theta}\right)d_{\downarrow}z
\eqno{(2.11)}
$$
on the phase space $\bold{C}\times\bold{C}\times H$ is
$\Gamma (2)\ltimes \bold{Z}^2$--invariant, has residue one
at its poles $z\equiv\dfrac{T_3}{2}\,\roman{mod}\,(1,\tau )$,
and therefore can be pushed down to the three--dimensional
space $F:=\Gamma (2)\ltimes \bold{Z}^2\setminus(\bold{C}\times\bold{C}\times H)$
fibered over the total space $E$ of the elliptic
pencil $\{ E_t\}$.
\endproclaim

We will use (2.11) first of all in order to identify $F$ with the phase
space described in the Introduction, 0.2. 
Here is the formal construction.

\smallskip

\proclaim{\quad 2.4. Lemma--Definition} Let $(\pi :E\to B,D_i)$
be an elliptic pencil with $\Gamma (2)$--rigidity, as in 0.2.
Then there exists an affine line bundle
$\lambda :F\to E$, and a relative 1--form
$\nu_F\in \Omega^1_{F/B}(\lambda^{-1}(D_3))$
such that the map of sheaves of affine lines over $\Cal{O}_E$
$\{$sections of $F$ over $E\} \to \Omega^1_{E/B}(D_3)$:
$s\mapsto s^*(\nu_F)$ identifies the sheaf of sections of
$F/E$ with that of forms with residue one $\Cal{F}\subset
\Omega^1_{E/B}(D_3)$. Moreover,

\smallskip

a). $(F=F(E,\pi ,\{D_i\}),\lambda ,\nu_F)$ is unique up to a
unique isomorphism over $E$.

\smallskip

b). $(\lambda : \bold{Z}^2\setminus (\bold{C}\times\bold{C}\times H)
\to \bold{Z}^2\setminus (\bold{C}\times H), 
\nu =\left( 2\pi iy+\dfrac{\theta_z}{\theta}\right)d_{\downarrow}z)$
is the $F$--space for the pencil $\{ E_{\tau}\}$ over $H$,
with $D_i\equiv \dfrac{T_i}{2}\,\roman{mod}\,(1,\tau ).$
\endproclaim

{\bf Proof.} Uniqueness follows from general nonsense. 
For existence, we give a standard \v{C}ech--type construction
which will be useful later.

\smallskip

Put $U_i=E\setminus D_i$ for $i=0,1,2.$

\smallskip

Localizing on $B$, we may and will assume that $\Omega^1_{E/B}$ is $\Cal{O}_E$--free.

\smallskip

Choose $\nu_i\in\Gamma (U_i,\Cal{F}), i=0,1,2$ (recall that
$\Cal{F}$ consists of relative 1--forms with residue 1 at $D_3$)
and take for $\nu_3$ a generator of $\Omega^1_{E/B}$ over $\Cal{O}_E.$

\smallskip

Define the alternating \v{C}ech 1--cocycle in
$Z^1((U_i),\Cal{O}_E)$ by
$$
f_{ij}=\frac{\nu_j -\nu_i}{\nu_3}\ \roman{on}\ U_i\cap U_j.
$$
Use it to glue together $U_i\times \bold{A}^1:$
$$
(x\in U_i\cap U_j,\, p_i\in\bold{A}^1)\mapsto
(x\in U_j\cap U_i,\, p_j=p_i+f_{ij}(x)).
$$
Denote by $F$ the resulting space, with projection $\lambda :
(x,p)\mapsto x$ on $E$. Denote by $\nu_F$ the form whose restriction
to $U_j\times \bold{A}^1$ is $\nu_j-p_j\nu_3.$
One easily checks the compatibility, so that $\nu_F$ is a
section of $\lambda^*(\Cal{F})\subset\Omega^1_{F/B}(\lambda^{-1}(D_3)).$

\smallskip

Clearly, $(F,\lambda ,\nu_F)$ satisfies the defining 
universal property. In fact, any section $\nu$ of $\Cal{F}$
on $U_i$ can be uniquely represented as a sum of $\nu_i$
and a unique regular differential, i.~e.
$\nu=\nu_i+f_i\nu_3,\ f_i\in\Gamma (U_i,\Cal{O}_E).$
Therefore $\nu$ is induced by $\nu_F$ on the section
locally given by
$U_i\to F: x\mapsto (x,f_i(x))$. 
We leave the last statement to the reader.

This finishes the proof.

\smallskip

Notice that $\lambda :F\to E$ has no global sections,
even over a single fibre of $E$, because there are no
differentials of the third kind with a single pole.
However, $F$ can be trivialized over $E\setminus D_i$
for any $i$ so that any Painlev\'e equation
with one nontrivial $\alpha_i$ effectively
lives on $E\times\bold{A}^1.$

\smallskip

Having thus described $(F,\nu_F )$, we can characterize
the whole space of Painlev\'e forms along the lines of 0.2.

\smallskip

\proclaim{\quad 2.5. Theorem} a). The form $\omega^{(0)}$
which in the $(y,z,\tau )$--coordinates is defined by
$$
\omega^{(0)} := 2\pi i\, (dy\wedge dz-ydy\wedge d\tau )
$$
is the unique closed holomorphic 2--form on $F$ satisfying
two conditions:

\smallskip

i). The restriction of $\omega^{(0)}$ to $\Cal{T}_{F/B}$
coincides with $d_{\downarrow}\nu_F$.

\smallskip

ii). The canonical lifts to $F$ of the multisections
of finite order of $E/B$ defined by
$$
z=e\tau +f,\ y=e;\ e,f\in \bold{Q}
\eqno{(2.12)}
$$
are leaves of of the null--foliation of $\omega^{(0)}$.

\medskip

b). The form $\omega_j$ on $E$ which in the $(z,\tau )$--coordinates
is defined by
$$
\omega_j:=\frac{1}{2\pi i}\wp_z(z+\frac{T_j}{2},\tau )dz\wedge d\tau
\eqno{(2.13)}
$$
is the unique closed meromorphic form on $E$ satisfying
two conditions:

\smallskip

iii). The divisor of $\omega_j$ is $\dfrac{D_kD_lD_m}{D_j^3},\ \{ j,k,l,m\}=
\{ 0,1,2,3\}$.

\smallskip

iv). If we identify $\Omega^2_E$ with $\pi^*(\Omega^1_{E/B})^{\otimes 3}$
with the help of the Kodaira--Spencer isomorphism
$d\tau \mapsto 4\pi i (d_{\downarrow}z)^2$, then in a formal
neighbourhood of $D_j$, $\omega_j$ becomes the cube of a
differential with constant residue $r$, where $r^3=-4.$
\endproclaim

{\bf Proof.} a). From (2.11) one sees that
$$
d_{\downarrow}\nu =2\pi i\,d_{\downarrow}y\wedge d_{\downarrow}z=
\left.\omega^{(0)}\right|_{\Cal{T}_{F/B}} .
$$
From (2.1) and (2.2) for $\alpha_i=0$ for $i=0,\dots ,3$ it follows that
(2.12) are solutions to this PVI.

\smallskip

Conversely, consider a holomorphic closed  2--form
$\tilde{\omega}^{(0)}$ enjoying properties i), ii).
Then
$$
\frac{1}{2\pi i}\tilde{\omega}^{(0)}=
dy\wedge dz+Edy\wedge d\tau +Gdz\wedge d\tau ,
$$
where $E, G$ are entire functions of $y,z,\tau $
with $E_z=-G_y.$ The respective equations of motion are
$$
\frac{dz}{d\tau}=-E(y,z,\tau ),\ \frac{dy}{d\tau}=G(y,z,\tau ).
\eqno{(2.14)}
$$
If (2.12) satisfy (2.14) for all $e,f\in\bold{Q}$,
we get that $E(e,e\tau +f,\tau)= -e$ for all
real $e,f$ by continuity, so that $E(y,z,\tau )\equiv -y$
by analyticity. Similarly, $G(e,e\tau +f,\tau )=0$
for all $e,f\in\bold{R}$ so that $G\equiv 0$, and
$\tilde{\omega}^{(0)}={\omega}^{(0)}.$

\smallskip

b). The divisor of (2.13) is well known. If the Kodaira--Spencer
isomorphism is normalized as above, $\omega_j$
becomes represented by the cubic differential
$$
2\wp_z(z+\frac{T_j}{2},\tau )(d_{\downarrow}z)^3=
\left( -\frac{4}{(z+\frac{T_j}{2})^3} +
O(z+\frac{T_j}{2}) \right)(d_{\downarrow}z)^3
$$
so that its formal cubic root near $z=-\frac{T_j}{2}$ exists 
and has residue $-\root 3\of{4}$. Any other cubic differential with
the same divisor can be obtained from ours by multiplication
by a function of $\tau$. Fixing the residue, we
lose this freedom.

\smallskip

\proclaim{\quad 2.6. Theorem} The Painlev\'e forms are exact.
More precisely,
$$
\omega (\alpha_0,\dots ,\alpha_3) =
d\Omega (\alpha_0,\dots ,\alpha_3)
\eqno{(2.15)}
$$
where $\omega (\alpha_0,\dots ,\alpha_3)$ is defined by (2.3), and
$$
\Omega (\alpha_0,\dots ,\alpha_3)=2\pi i\,(ydz-\frac{1}{2}y^2d\tau)
+d\roman{log}\,\theta(z,\tau )
+2\pi i\,G_2(\tau )d\tau +
$$
$$
+\frac{1}{2\pi i}\sum_{j=0}^3\alpha_j\wp (z+\frac{T_j}{2},\tau )d\tau
\eqno{(2.16)}
$$
is a $\Gamma (2)\ltimes\bold{Z}^2$--invariant meromorphic 1-form
with poles of the second order at $D_j$.

Here
$$
G_2(\tau ):=-\frac{1}{24}+\sum_{n=1}^{\infty}
(\sum_{d/n}d)e^{2\pi in\tau} .
\eqno{(2.17)}
$$
\endproclaim 
\smallskip

{\bf Proof.} Only $\Gamma (2)\ltimes\bold{Z}^2$--invariance needs to
be checked.
This is a straightforward calculation using (2.4), (2.5), (2.7),
(2.8), and the pseudo--modular property of $G_2(\tau )$:
$$
G_2\left(\frac{a\tau +b}{c\tau  +d}\right)=
(c\tau +d)^2G_2(\tau )-\frac{c(c\tau +d)}{4\pi i} .
$$
We leave it to the reader.

\smallskip

\proclaim{\quad 2.7. Theorem} On a PVI phase space $(F,\lambda ,\nu_F)$,
denote by $D$ the divisor of the zeroes of $\nu_F$ considered
as a section of the invertible sheaf $\lambda^*(\Omega_{E/B}^1(D_3)).$
Then:

\smallskip

a). $D$ is a section of $\lambda : F\setminus\lambda^{-1}(D_3)\to
E\setminus D_3.$

\smallskip

b). In the space of Painlev\'e 2--forms, there exists a unique
form identically vanishing on $D$. It corresponds to the
point $(\alpha_0,\dots ,\alpha_3)=(0,0,0,\frac{1}{2})$ (which
is the $\bold{P}^2$--point up to a renumbering and a
Landin transform, cf. 0.1 and 1.6.)

\smallskip

c). $D$ is generically transversal to the null--leaves of any
PVI except for $(0,0,0,\frac{1}{2})$, and so can serve as
a common space of initial conditions for these equations
(cf. [O3] for a framework for the more precise analysis.)
\endproclaim

\smallskip

{\bf Proof.} a). From the construction given in the proof of Lemma 2.4
one sees that the equation of $D$ in $U_j\times\bold{A}^1$
is $p_j=\nu_j/\nu_3.$ Since $\nu_3$ is everywhere
invertible, and the only pole of $\nu_j$ is $D_3$, $D$ is
a section of $\lambda$ outside $D_3$.

\smallskip

b). Since the difference of any two forms in the Painlev\'e
space on $F$ is lifted from $E$, its restriction to $D$
can vanish identically only if these forms coincide.
Hence at most one form can vanish on $D$ identically. To exhibit the
one corresponding to $(\alpha_0,\dots ,\alpha_3)=
(0,0,0,\frac{1}{2})$ we will prove a slightly stronger statement,
that it is a differential of a form vanishing on $D.$
Put
$$
\Omega_0:=2\pi i\,[ydz-\frac{1}{2}y^2d\tau]+d\roman{log}\,
\theta(z,\tau )+
\frac{i}{4\pi}\frac{\partial^2}{\partial z^2}\,\roman{log}\,
\theta (z,\tau )d\tau .
\eqno{(2.18)}
$$
Then we can consecutively check that it is $\Gamma (2)\ltimes
\bold{Z}^2$--invariant and that $d\Omega_0=\omega (0,0,0,\frac{1}{2}).$
For the latter, use the identity
$$
\frac{\partial^2}{\partial z^2}\,\roman{log}\,
\theta (z,\tau )=
-\wp (z+\frac{1+\tau}{2},\tau )+\varphi (\tau )
$$
where the precise form of $\varphi (\tau )$ is inessential here.

\smallskip

On the other hand, from the heat equation
$$
\theta_{\tau}(z,\tau )=\frac{1}{4\pi i}\,\theta_{zz}(z,\tau )
$$
it follows that
$$
\Omega_0=\left[2\pi i\,y+\frac{\theta_z}{\theta}\right]dz
-\frac{1}{2}\left[2\pi i\,y^2+\frac{i}{2\pi}
\left(\frac{\theta_z}{\theta}\right)^2\right]d\tau =
$$
$$
=\left( 2\pi i\,y+\frac{\theta_z}{\theta}\right)
\left[dz-\frac{1}{4\pi i}\left(2\pi i\,y-\frac{\theta_z}{\theta}
\right)d\tau \right] .
\eqno{(2.19})
$$
Comparing this with (2.11), one sees that $\Omega_0$ vanishes on $D.$

\smallskip

c). From the previous discussion, it follows that if 
$\omega\ne d\Omega_0$, then the restriction of 
$\omega$ to $D$ is generically of rank 2, so that
its null--foliation is generically transversal to $D$.

\medskip

{\bf 2.8. The structure of the phase space in algebraic coordinates.}
In this subsection, we will work out the basic formulas
on the algebraic model (0.3).

\smallskip

{\bf 2.8.1. The vertical coordinate.} According to the proof
of Lemma 2.4, the natural vertical (over $E$) coordinate
on the algebraic model of $F\setminus \lambda^{-1}(D_3)$ is
$$
U:=\frac{\nu}{d_{\downarrow}X/Y} .   \eqno{(2.20)}
$$
From (2.11) and (1.14) one finds its expression through
elliptic functions:
$$
U=\roman{(push\ down\ of)}\ 
\frac{ 2\pi i\,y+\dfrac{\theta_z}{\theta}}{2(e_2(\tau )-e_1(\tau ))^{1/2}}.
\eqno{(2.21)}
$$
In particular, the equation of $D$ is simply $U=0.$

\smallskip

We will now identify the Painlev\'e forms, using the classical
parameters $(\alpha ,\beta ,\gamma ,\delta )$ rather than $\alpha_i$.

\smallskip

\proclaim{\quad 2.8.2. Theorem} We have
$$
\Omega (\alpha ,\beta ,\gamma ,\delta )=
U\frac{dX}{Y}-U^2\frac{dt}{t(t-1)}+
$$
$$
+\frac{1}{2t(t-1)}
\left( \alpha X-\beta \frac{t}{X} - \gamma \frac{t-1}{X-1}
-\delta \,\frac{t(t-1)}{X-t}\right) dt,
\eqno{(2.22)}
$$
$$
\omega (\alpha ,\beta ,\gamma ,\delta )=
dU\wedge\frac{dX}{Y}-\frac{U}{2(X-t)Y}\,dX\wedge dt
-2UdU\wedge\frac{dt}{t(t-1)}+
$$
$$
+\frac{1}{2t(t-1)}
\left( \alpha +\beta \frac{t}{X^2} + \gamma \frac{t-1}{(X-1)^2}
+\delta \frac{t(t-1)}{(X-t)^2}\right) dX\wedge dt .
\eqno{(2.23)}
$$
\endproclaim
\smallskip

{\bf Proof.} The main task is to show that (2.22) holds
for $\alpha =\beta =\gamma =\delta =0.$ In fact, the part involving
$\alpha ,\beta ,\gamma ,\delta $ can be treated in the same way as
at the end of the proof of Theorem 1.4, and (2.23) is
then
obtained by derivation.

\smallskip

Now, the form $\Omega$ corresponding to $\alpha =\beta =\gamma =\delta =0$  
is precisely $\Omega_0$ from (2.19). Thus, we have to prove
that
$$
\Omega_0=U\frac{dX}{Y}-U^2\frac{dt}{t(t-1)} .
\eqno{(2.24)} 
$$
Using (2.19) and (2.21), we find
$$
\Omega_0=\frac{2\pi i\,y+\dfrac{\theta_z}{\theta}}{2(e_2-e_1)^{1/2}}
\cdot 2(e_2-e_1)^{1/2}\left[ \frac{i}{4\pi}\,
\frac{2\pi i\,y+\dfrac{\theta_z}{\theta}}{2(e_2-e_1)^{1/2}}
\cdot 2(e_2-e_1)^{1/2}d\tau
+dz +\frac{1}{2\pi i}\,\frac{\theta_z}{\theta}d\tau  \right]=
$$
$$
=U^2\,\frac{i}{\pi}(e_2-e_1)d\tau +
U\cdot 2(e_2-e_1)^{1/2}\left[ dz +\frac{1}{2\pi i}\,\frac{\theta_z}{\theta}
d\tau  \right] .
\eqno{(2.25)}
$$
From (1.13) one can deduce that
$$
\frac{i}{\pi}(e_2-e_1)d\tau =-\frac{dt}{t(t-1)} .
$$
Comparing (2.25) and (2.24), one sees that it remains to prove that
$$
\frac{dX}{Y}=2(e_2-e_1)^{1/2}\left[ dz +\frac{1}{2\pi i}\,\frac{\theta_z}{\theta}
d\tau  \right] .
\eqno{(2.26)}
$$
Now, from (1.13) we obtain
$$
\frac{dX}{Y}=d\left(\frac{\wp -e_1}{e_2-e_1}\right)\cdot
\frac{2(e_2-e_1)^{3/2}}{\wp_z}=
$$
$$
=2(e_2-e_1)^{1/2}dz-2\frac{\wp_{\tau}-e_{1\tau}}{(e_2-e_1)^{1/2}\wp_z}d\tau .
\eqno{(2.27)}
$$
Taking the difference of the right hand sides of (2.26) and (2.27),
we first check that it cannot depend on $z$, because
a calculation shows that $\dsize d\left(\frac{dX}{Y}-\mu\right)=0$,
where
we temporarily denoted by $\mu$ the right hand side of (2.26).
Put now $\dsize\frac{dX}{Y}-\mu=\varphi (\tau )d\tau .$
Then we can calculate $\varphi (\tau )$ by restricting this identity
to the divisor $D_1: X=0$ or equivalently, $z=1/2.$
We get
$$
\varphi (\tau )=\frac{1}{2\pi i}\,\frac{\theta_z(1/2,\tau )}{\theta (1/2,\tau )}
=0
$$
finishing the proof.

\medskip

{\bf 2.8.3. PVI in the $(U,X,Y,t)$--space and the canonical
lifts of the multisections of the finite order.}
From (2.23) one deduces the following equations of motion:

\smallskip

$\dsize\frac{dX}{dt}=\frac{2UY}{t(t-1)},$
$$
\frac{dU}{dt}=-\frac{U}{2(X-t)}
+\frac{Y}{2t(t-1)}\left( \alpha +\beta \frac{t}{X^2} + \gamma \frac{t-1}{(X-1)^2}
+\delta \frac{t(t-1)}{(X-t)^2}\right) .
\eqno{(2.28)}
$$
In particular, if $(X(t),Y(t))$ is a multisection of finite order, hence
a solution of the $(\alpha =\beta =\gamma =0,\, \delta =1/2)$--equation,
then from the first equation (2.28) we see that its lift to $F$ is given by
$$
U(t)=\frac{t(t-1)}{2}\,\frac{X^{\prime}(t)}{Y(t)}.
\eqno({2.29)}
$$

\vskip 1cm

\centerline{\bf \S 3. Symmetries and special solutions}

\bigskip

{\bf 3.1. Reduced phase space and enhanced moduli space.}
The discrete symmetries of PVI of infinite order were discovered in the context
of isomonodromic deformations by Schlessinger and rediscovered
many times afterwards (cf. [JM].) We review here 
the Okamoto's treatment [O2] which nicely fits in our framework.

\smallskip

Our phase space $F$ has an obvious $\bold{Z}_2$-symmetry
induced by the inversion map on the fibers of $E$:
$$
(y,z,\tau )\mapsto (-y,-z,\tau ),\ (U,X,Y,t)\mapsto (-U,-Y,X,t).
$$
Each Painlev\'e form and the respective equations of motion are invariant
w.r.t. this symmetry. We delete eventual poles and consider the reduced
phase space $\dsize F_0:=(F\setminus\cup_{i=0}^3 D_i)/\bold{Z}_2$. In this
section we will work with  the algebraic $(U,Y,X,t)$--model. Then $F_0$
has an obvious structure of affine algebraic variety.

\smallskip

We also replace the moduli space $P_0:=\roman{Spec}\,\bold{C}[\alpha_0,\dots ,
\alpha_3]$ by its cover
$$
P:=\roman{Spec}\,\bold{C}[a_0,\dots ,a_3],\ a_i^2=2\alpha_i.
$$
Finally, we introduce the pair
$(\Phi :=F_0\times P,\omega)$, where $\omega$ is the (relative over $P$)
closed regular algebraic 2--form on $\Phi$ denoted 
$\omega (\alpha ,\beta ,\gamma ,\delta )$ in (2.23).
This pair is an algebraic model of the space of all PVI equations.

\medskip

{\bf 3.2. Symmetries.} Denote by $W$ the group of symmetries of $P$
generated by the following transformations:

\smallskip

a). $(a_i)\mapsto (\varepsilon_ia_i)$, where $\varepsilon_i=\pm 1.$

\smallskip

b). Permutaions of $(a_i).$

\smallskip

c). $(a_i)\mapsto (a_i+n_i)$, where $n_i\in\bold{Z}$
and $\dsize \sum_{i=0}^3n_i\equiv 0\,(2)$.

\smallskip

\proclaim{\quad 3.2.1. Theorem (Okamoto [O2])} All transformations
in $W$ can be lifted to the birational transformations
of $\Phi =F_0\times P$ preserving the equations of motion
defined by $\omega$.
\endproclaim

Sign changes of $a_i$ can be extended by identity on $F_0$.
The action of $S_4$ on $E$ was described in 1.6.
To lift it to $F$, it suffices to remark that the four affine sheaves
of differentials with a single pole and residue 1 at $D_i$
can be pairwise identified by adding
$\frac{1}{2}d\,\roman{log}\,f_{ij},\ \roman{div}\,f_{ij}=
D_j^2/D_i^2.$  

\smallskip

The whole
group $W$ is generated by these elements and one shift
$(a_i)\mapsto (a_i+\delta_{i0}+\delta_{i3})$, hence it suffices
to construct its lifting. I will briefly sketch Okamoto's
ingenious argument for doing this.

\smallskip

I start with comparing notation. Okamoto's $q,t$ are ours $X,t$.
The vertical coordinate in the phase space which
Okamoto denotes $p$ can be identified as
$$
p=\frac{U}{Y}+\frac{1}{2}\left(\frac{a_1}{X}+\frac{a_2}{X-1}+\frac{a_3-1}{X-t}
\right).
\eqno{(3.1)}
$$
The verification reduces to a somewhat tedious calculation,
showing that (3.1) transforms Okamoto's equations of motion
([O2], (1.5), p. 349) into ours (2.28). It is useful to
remember that Okamoto's parameters $(\kappa_{\infty},
\kappa_0,\kappa_1,\theta )$ are ours $(a_0,a_1,a_2,a_3).$

\smallskip

Okamoto introduces an auxiliary function $h$ ((1.6), p. 349),
which in our coordinates is
$$
h=U^2+\frac{1}{4}\left[-a_0^2X-a_1^2\frac{t}{X}
+a_2^2\frac{t-1}{X-1} -(a_3-1)^2\frac{t(t-1)}{X-t} \right]-
$$
$$
-\frac{1}{4}(a_3-1)^2t+ \frac{1}{8}[a_0^2+a_1^2
-a_2^2 +(a_3-1)^2].
\eqno{(3.2)}
$$
We need also the Painlev\'e flow on $\Phi$ given by the 
total time derivative
$$
\Cal{D}:=\partial_t+\frac{2UY}{t(t-1)}\partial_X-
$$
$$
-\left[
\frac{U}{2(X-t)}
-\frac{Y}{4t(t-1)}\left (a_0^2 -a_1^2 \frac{t}{X^2} +a_2^2 \frac{t-1}{(X-1)^2}
-(a_3^2-1) \frac{t(t-1)}{(X-t)^2}\right)\right]\partial_U.
\eqno{(3.3)}
$$ 
This is a restatement of (2.28).

\smallskip

Now Okamoto's description of the shift can be summarized
as follows.

\smallskip

{\it The action of the shift upon $h$ given explicitly by
$$
h\mapsto h-X(X-1)\left( \frac{U}{Y}+\frac{a_1}{2X}+\frac{a_2}{2(X-1)}+\frac{a_3-1}{2(X-t)}\right) +
$$
$$
+\frac{1}{2}(-a_0+a_1+a_2+a_3-1)X+\frac{1}{4}(a_0-2a_1-a_3+1)
$$
has a unique birational extension to the whole affine ring of $\Phi$
compatible with $\Cal{D}$. }

\smallskip

The proof given by Okamoto is a calculation. He shows that the ring
homomorphism 
$$
\bold{C}[a_0,\dots ,a_3;h,k,l]\to \roman{affine\ ring\ of\ }\Phi
$$
defined by $h\mapsto h,\, k\mapsto \Cal{D}h,\, l\mapsto \Cal{D}^2h$,
after a localization becomes surjective. Its kernel is generated by an explicit
polynomial relation (see [O2], p. 349, Prop. 1.1.)
The symmetries of this polynomial relation are slightly more
visible than those of the initial setting.

\smallskip

The geometric meaning of this proof in the context of elliptic pencils remains unclear to me. On the level of the complete phase space $F\times P$,
Okamoto's map becomes a correspondence.

\medskip

{\bf 3.3. Special solutions.} The points of the Painlev\'e moduli
space can be roughly divided into four groups, according to
the dimension of the space of solutions reducible to
``classical'' functions. This is a purely experimental 
classification, since to the author's knowledge, no precise definition of this notion led to a precise classification picture.
Nevertheless, it seems worthwhile to summarize a part of what is known.

\smallskip

Generally, some classical solutions  at a point of $P$ are constructed  
directly. Afterwards new solutions can be generated in principle
 by applying transformations from $W$ and Landin's transform
(which in the $(a_i)$ coordinates is
$(a_0,a_1,a_0,a_1)\leftrightarrow (2a_0,2a_1,0,0).)$ Especially interesting
are {\it algebraic} solutions: those for which $X(t)$ is an algebraic
function of $t$. Symmetries (including Landin) preserve
the algebricity.

\smallskip

a). {\it Equations with classical general solutions.}
The basic point for them is the null--point $a_i=0$ for all $i.$
In the $(z,\tau )$--space we get simply $z= e\tau +f,\,e,f\in\bold{C}.$
Algebraic solutions are obtained precisely for $e,f\in\bold{Q}.$
They are rigid (not deformable), but in a certain sense dense
in the set of all solutions.

\smallskip

Applying shifts from $W$, we get infinitely many
classically completely solvable equations:
$(a_i)=(n_i),\ (\alpha_i)=(\dfrac{n_i^2}{2})$, where $n_i$
are integral and the sum of $n_i$ is even. Let $L$
be the lattice of such vectors. Inverse
Landin transform applied to $(1,1,0,0)$ then shows that the points $(a_i)$ in
$(\frac{1}{2},\frac{1}{2},\frac{1}{2},\frac{1}{2})+L$ are also classically solvable;
this includes Hitchin's equation.

\smallskip

I do not know of other PVI equations with this property.

\smallskip

b). {\it Equations with one--dimensional families of classical
solutions.} The basic point here is the $\bold{P}^2$--point
$(a_i)=(0,0,0,1).$ One family of solutions is obvious
in the algebraic coordinates (2.28): $X=\roman{const}, U=0.$
These solutions have a clear geometric meaning in our phase
space $F$: they form the foliation of the
divisor $D$ formed by the null--leaves entirely contained in $D$:
cf. Theorem 2.7.

\smallskip

It is interesting that this time algebraic solutions are not rigid.
If they look somewhat plain, they become more sophisticated
on other elements of the orbit $(0,0,0,1)+L$ and on the
Landin transforms $(0,0,\frac{1}{2},\frac{1}{2})+L$ and $(\frac{1}{4},\frac{1}{4},\frac{1}{4},\frac{1}{4})+L.$

\smallskip

More interesting one--dimensional family of solutions 
can be constructed for any $(a_i)$ belonging to the
hyperplane $a_0+a_1+a_2+a_3=1$. They are expressed through
Gauss hypergeometric equations: see [O2],
p. 373--374. Again, $W$ and Landin
generate infinitely many new families.

\smallskip

c). Hitchin [H1] and Dubrovin [D] constructed isolated algebraic solutions
using respectively twistor geometry and Frobenius manifolds.

\smallskip

Our last remark concerns some similarity between
the (generalized) Lam\'e potentials in the theory of KdV--type
equations and our classically integrable
potentials of the non--linear equation (2.2).
According to [TV], the former are of the form
$$
\sum_{j=0}^3\frac{n_j(n_j+1)}{2}\wp(z+\frac{T_j}{2},\tau ),
$$
whereas according to our discussion
the latter have coefficients (proportional to)
$(n_j^2)/2$ or $(n_j+\frac{1}{2})^2/2$. Is there a direct connection
between the two phenomena?

\bigskip

\centerline{\bf References}

\bigskip

[D] B.~Dubrovin. {\it Geometry of 2D topological field theories.}
In: Springer LNM, 1620 (1996), 120--348.

\smallskip

[DFI] P.~Di Francesco, C.~Itzykson. {\it Quantum intersection
rings.} In: The Moduli Space of Curves, ed. by R.~Dijkgraaf,
C.~Faber, G.~van der Geer, Progress in Math., vol. 129,
Birkh\"auser 1995, 149--163.

\smallskip

[F] R.~Fuchs. {\it \"Uber lineare homogene Differentialgleichungen
zweiter Ordnung mit im endlich gelegene wesentlich
singul\"aren Stellen.} Math. Ann., 63 (1907), 301--321.

\smallskip

[G] B.~Gambier. {\it Sur les \'equations diff\'erentielles
du second ordre et du pr\'emier degr\'e dont l'int\'egrale
g\'en\'erale est \'a points critiques fixes.} CR Ac. Sci. Paris,
142 (1906), 266--269.

\smallskip

[H1] N.~Hitchin. {\it Poncelet polygons and the Painlev\'e equations.}
In: Geometry and Analysis, ed. by S. Ramanan,
Oxford University Press, Bombay, 1995.

\smallskip

[H2] N.~Hitchin. {\it Twistor spaces, Einstein metrics and
isomonodromic
deformations.} J. Diff. Geom., 3 (1995), 52--134.

\smallskip

[H3] N.~Hitchin. {\it Frobenius manifolds (notes by D. Calderbank.)}
Preprint, 1996.

\smallskip

[JM] M.~Jimbo, T.~Miwa. {\it Monodromy preserving deformation
of linear ordinary differential equations with rational
coefficients II.} Physica 2D (1981), 407--448.

\smallskip

[M] Yu.~Manin. {\it Rational points of algebraic curves over
functional fields.} AMS Translations, ser. 2, vol. 50 (1966),
189--234.

\smallskip

[O1] K.~Okamoto. {\it Isomonodromic deformation and
Painlev\'e equations, and the Garnier system.}
J. Fac. Sci. Univ. Tokyo, Sect. IA Math., 33 (1986), 575--618.

\smallskip

[O2] K.~Okamoto. {\it Studies in the Painlev\'e equations I.
Sixth Painlev\'e equation PVI.} Annali Mat. Pura Appl.,
146 (1987), 337--381.

\smallskip

[O3] K.~Okamoto. {\it Sur les feuilletages associ\'es aux
\'equation du second ordre \`a points critiques fixes
de P. Painlev\'e. Espaces de conditions initiales.}
Japan J. Math., 5:1 (1979), 1--79.

\smallskip

[TV] A.~Treibich, J.--L.~Verdier. {\it Rev\^etements
tangentiels et sommes de 4 nombres triangulaires.}
C.R. Ac. Sci. Paris, s\'er. I Math., 311 (1990), 51--54.

\enddocument